\documentclass{ifacconf}

\usepackage{graphicx}      
\usepackage{natbib}        
\usepackage{tikz}
\usepackage{subfigure}
\newtheorem{definition}{Definition}
\def\Ii{{\bf 1}} 
\def\PP{{\bf P}} 
\def\EE{{\bf E}} 

\begin{document}
\begin{frontmatter}
\title{Extremes Control of Complex Systems With Applications to Social Networks}

\thanks[footnoteinfo]{The 
author  was partly supported by the Russian Foundation for Basic Research, grant 13-08-00744 A.}

\author[Markovich]{Natalia M. Markovich}
\address[Markovich]{Institute of Control Sciences, Russian Academy of Sciences, Moscow, Russia (e-mail: markovic@ipu.rssi.ru)}
\begin{abstract}
The control and risk assessment in complex information systems require to take into account extremes arising from nodes with large node degrees. Various sampling techniques like a Page Rank random walk, a Metropolis-Hastings Markov chain and others serve to collect information about the nodes.
The paper contributes to the comparison of sampling techniques in complex networks by means of the first hitting time, that is the minimal time required to reach a  large node. Both the mean and the distribution of the first hitting time is shown to be determined by the so called extremal index. The latter indicates a dependence measure of extremes and also reflects the cluster structure of the network. The clustering is caused by dependence between nodes and heavy-tailed distributions of their degrees. Based on extreme value theory  we estimate the mean and the distribution of the first hitting time and the distribution of node degrees by  real data from social networks. We demonstrate the heaviness of the tails of these data using appropriate tools. The same methodology can be applied to other complex networks like peer-to-peer telecommunication systems.
 \end{abstract}
\end{frontmatter}
\section{Introduction}
\begin{figure}[htb]
\centering
\includegraphics[width=0.46\textwidth]{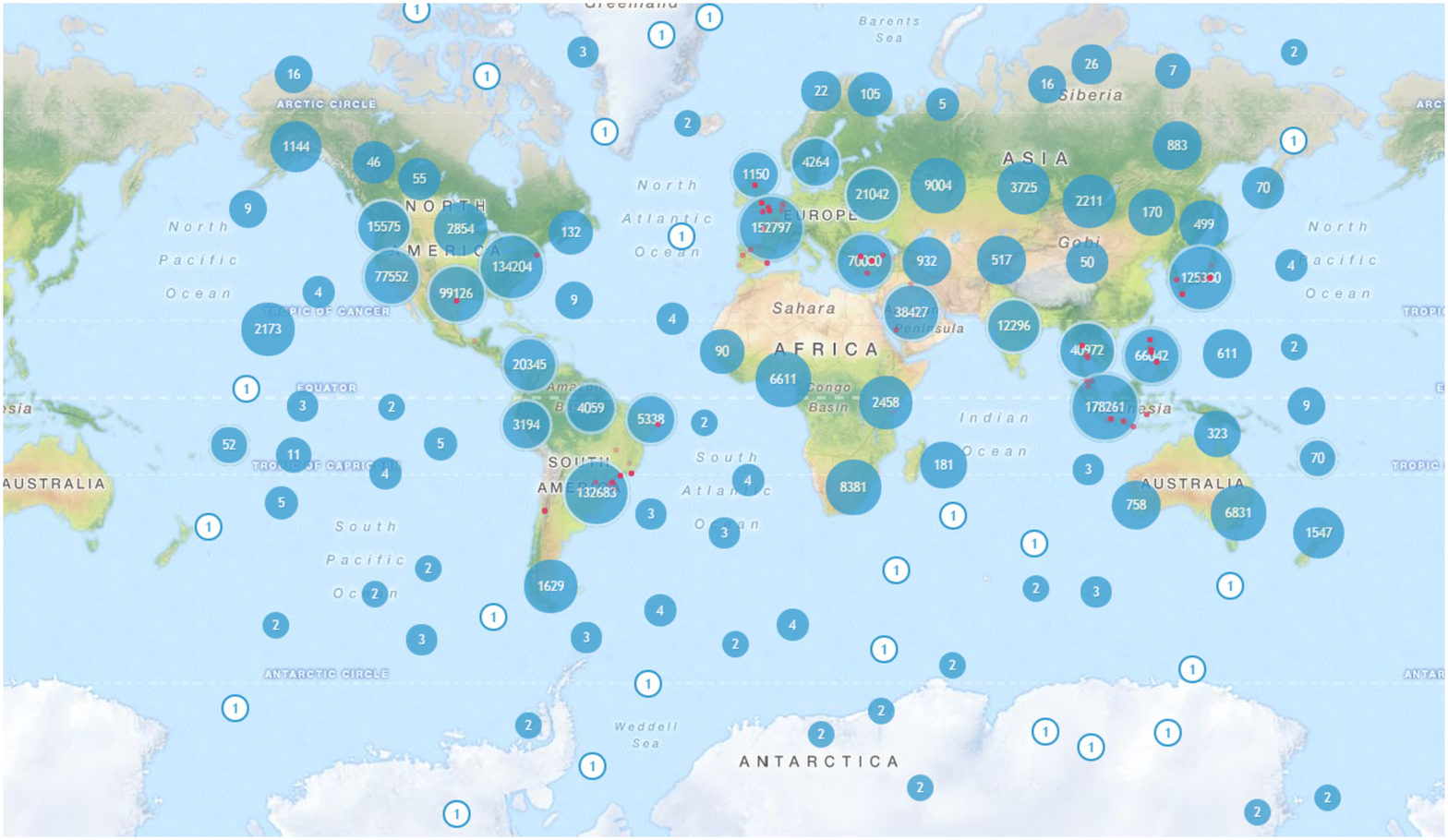}
\caption{Clusters of nodes on a Twitter map \cite{twit}.}
\label{fig:1}
\end{figure}

Modern complex networks like online social (OSN), peer-to-peer (P2P) and content-centric networks and the world wide web (WWW)   are in general nonlinear information systems. The control and risk assessment in complex information systems require to take into account extremes arising from nodes with large node degrees. The giant extremal nodes impact on the work and the development of the whole system more than small nodes.
The investigation of all nodes is costly since the networks are very large. Thus, sampling techniques via crawling  are proposed as tools to collect node samples. Uniform and random walk sampling, PageRank (\cite{Avrachenkov}), non-backtracking random walk with re-weighting (NBRW) (\cite{ChulHo}), the random walk Metropolis and Metropolis-Hastings  algorithms (\cite{Metropolis}, \cite{Hastings})  give examples of possible approaches.
\\
The giant nodes surrounded by smaller nodes build clusters of connected nodes. Within the clusters, node degrees may exceed sufficiently large thresholds $\{u_n\}$. By the extreme value theory the node degrees exceeding $u_n$ such that $u_n\to\infty$ as sample size $n\to\infty$ build a compound Poisson process, \cite{Beirlant}, \cite{Leadbetter83}. Roughly speaking, tops of the clusters become independent and determine independent clusters of the network. It is visible for a Twitter network taken as an example in Fig. \ref{fig:1}.
With this respect, it is important to evaluate the first hitting time, that is the minimal time required to reach a  large node. This allows us to disseminate information and advertisement more effectively and to upload it directly to top-nodes of such clusters.
\\
The extremal index is a key characteristic of cluster extremes. Its reciprocal approximates the mean cluster size, i.e. the number of exceedances of the threshold per cluster, \cite{Leadbetter83}. It determines the first hitting time and its distribution and mean, \cite{Roberts}, \cite{Markovich}. The extremal index allows us to represent the distribution of the maximal node degree (\ref{6}) and its quantiles. Thus, the estimation of the extremal index is one of the subjects of the paper.
\\
Another problem that is related to the extremal index estimation is given by the detection of the heaviness of tails of the node degree distribution. The presence of  heavy tails may dramatically impact on the first hitting time of the sampling technique and its effectiveness. The power law distribution which has asymptotically  a Pareto tail  is widely applied to model the node degree distribution, \cite{Litvak}, \cite{Newman}. However, such models may fit  the  distributions unsatisfactory and do not satisfy nonparametric tests, \cite{Litvak}. 
the reason is that distributions may include mixtures of heavy- and light-tailed distributions. This applies also to the popularity of a content transmitted through content-centric networks, \cite{Rossi}.
\\
We focus here on the Metropolis algorithm  which constructs nonlinear time-reversible Markov chains with a given, desired stationary distribution $\pi(x)$, \cite{Andrieu}.  For a given graph of the network there are potentially many irreducible\footnote{This means that every node in the network is reachable in a finite time with a positive probability.} Markov chains (or random walks) preserving the same stationary distribution $\pi$.
To select the best one,  criterion like the mixing time (\cite{ChulHo}) and related to it the second largest eigenvalue  of the
associated random walk transition probability matrix, \cite{Avrachenkov} as well as the convergence rate of a Markov chain (\cite{MengersenTweedie}) are usually applied. \\
On this respect, heavy-tailed distributions $\pi(x)$ of a Metropolis random walk generate specific problems. The convergence rate of such Markov chain is not a geometric but a polynomial one, \cite{RobertsSmith}. This leads to an infinitely long first hitting time to reach a node with a large degree. The latter time is finite in case of a light-tailed $\pi(x)$.
\\
The objectives of the paper are the following. Summarizing theoretical achievements regarding extremes of stochastic sequences, we apply them to real data sets of social networks. We make conclusions regarding the possible effectiveness of
Markov chains generated by  the Metropolis algorithm with a heavy-tailed stationary distribution $\pi(x)$.
\\
To this end
(1) we evaluate the node degree distribution by real data and test it regarding the presence of heavy tails by appropriate tools, 
(2)
we investigate extreme value statistics such as the extremal index  and the first hitting time, namely, its  mean and distribution  by real data of social networks.
\\
The paper is organized as follows. In Section \ref{Defin}  definitions and related theoretical results concerning  the extremal index and the first hitting time are given. The power law model is given in Section \ref{Power Law}. The  random walk Metropolis algorithm is described in Section  \ref{Metropolis}.
Applications to real data sets of two social networks and the estimation of the extremal index and the first hitting time are given in Section \ref{Model}.  Conclusions are stated in Section \ref{Conclus}.
\section{Definitions and related work}
\subsection{The extremal index and the first hitting time}\label{Defin}
Let $\{X_n\}_{n\ge 1}$ be a stationary sequence with marginal distribution function $F(x)$ and $M_n=\max\{X_1,...,X_n\}$. One may involve that $\{X_n\}$ is a sequence of node degrees.
\begin{definition}\cite{Leadbetter83} The stationary sequence  $\{X_n\}_{n\ge 1}$ is said to have the extremal index
$\theta\in[0,1]$ if
for each constant $0<\tau <\infty$ there is a sequence of real numbers (thresholds) $u_n=u_n(\tau)$ such that
\begin{equation}\label{1}\lim_{n\to\infty}n(1-F(u_n))=\tau \qquad\mbox{and}\end{equation}
\begin{equation}\label{2}\lim_{n\to\infty}P\{M_n\le u_n\}=e^{-\tau\theta}\end{equation}
hold.
\end{definition}
 $\theta$ is a dependence measure of extremes in the following sense. For a sufficiently large sample size $n$ and the threshold sequence  $\{u_n\}$
\begin{equation}\label{6}\PP\{M_n\le u_n\}=F^{n\theta}(u_n)+o(1)\end{equation}holds. If $\{X_n\}$ are independent random variables then
  $\theta=1$  as far as  $\theta\approx 0$ corresponds to a strong dependence.
    If  $\pi(x)$ behaves as a power function, i.e. if
    \begin{equation}\label{6a}\pi(x)\sim cx^{-\alpha}\end{equation}
    for some $c>0$ and $\alpha>0$, the Metropolis algorithm gives an  example of the pathological case $\theta= 0$ of total dependence, \cite{Roberts}.
         \\
     The  extremal index $\theta$ of  $\{X_n\}$ indicates the first hitting time  $T_n$  to exceed level $u_n$, \cite{Roberts}.
     \begin{definition} The first hitting time $T^*=T^*(u_n)$ of the threshold $u_n$ is determined by the following expression
     \[P\{T^*(u_n)=j+1\}= P\{M_j\le u_n, X_{j+1}>u_n\}, \] $j=0,1,2,...$, $M_{0}=-\infty$, \cite{Markovich}.
     \end{definition}
     Since
 $u_n$ is selected according to (\ref{1}) it follows that $\PP\{X_n > u_n\}$ is asymptotically equivalent to $1/n$. Notice, that it holds
 \[\PP\{M_k\le u_n\}=\PP\{T^*>k\}.\]
  Hence, we get\footnote{The $'$$\sim$$'$ means asymptotically equal to or $f (x)\sim g(x)$ $\Leftrightarrow$ $f (x)/g(x) \to 1$ as $x\to a$, $x\in M$.} \[\PP\{T^*/n>k/n\}\sim e^{-\theta k\PP\{X_n > u_n\}}\sim e^{-\theta k/n}\] and it follows
 \[\lim_{n\to\infty}\PP(T^*/n > x) = e^{-\theta x}\] for positive $x$. It follows
\begin{equation}\label{10}\lim_{n\to\infty}E(T^*/n)=1/\theta.\end{equation}
 $E(T^*/n)$ determines the mean  first hitting time to find a node with the degree exceeding a sufficiently large level.
 (\ref{10}) implies, that  the smaller $\theta$, the longer it takes to reach a node with a large degree.
 \\
    The latter result is specified in \cite{Markovich}. More exactly, it holds
    \begin{equation}\label{9}\lim_{n\to\infty}\rho_nET^*(x_{\rho_n})=1/\theta^3,\end{equation}
    where the quantile $x_{\rho_n}$ of the level $1-\rho_n$ of the underlying sequence\footnote{This implies that $\PP\{X_1>x_{\rho_n}\}=\rho_n$ holds.} $\{X_n\}$ is taken as the threshold $u_n$. Since $\rho_n\sim \tau/n$ according to (\ref{1}) the result (\ref{9}) does not contradict (\ref{10}).
 \\
 For example, if $\theta=1/2$ then it will take eight times  longer  to find the node degree $D_i=u_n$ than to arrive at extreme levels by an independent sequence.
    \\
    The normalized distribution of $T^*(x_{\rho_n})$ is derived to be geometric with the probability equal to $\rho_n\theta$, i.e.
     \begin{equation}\label{11}P\{T^*(x_{\rho_n})=j\}\sim  \frac{\rho_n\theta}{\theta^2}(1-\rho_n\theta)^{j-1}\end{equation} holds under a specific mixing condition, \cite{Markovich}.
\subsection{Power law}\label{Power Law}
 We consider node degrees of real social networks
  corresponding to indirected graphs. In this case, in- and out- degrees coincide. The node degree distributions are believed to follows power laws (\ref{6a})
  and for Web $\alpha=1.1$ for in-degree and PageRank, and $\alpha\approx 2$ for out-degree, \cite{VolkovichLitvakZwart}, \cite{volkovich2010}.  
  \\
  Despite the node degree is a discrete random variable, it is a common approach to consider a simpler continuous Pareto analogue of its distribution, \cite{Clauset}.
  The latter belongs to the class of heavy-tailed regularly varying distributions with the tail function
\begin{equation}\label{3}\PP\{X>x\}=\ell(x)x^{-\alpha},\quad x>0,\end{equation}
where $\alpha$ denotes the tail index responsible for the heaviness of the tail, and $\ell(x)$ is a slowly varying function, that is, for $x>0$, $\ell(tx)/\ell(t)\to 1$ as $t\to\infty$. In practice, the latter model  may fit the tail of the distribution  but unlikely the body of the distribution. Usually, we do not know $\ell(x)$ of the appropriate model. It could be any positive constant or logarithm. The model (\ref{3}) is sensitive to the estimation of the tail index $\alpha$, \cite{VolkovichLitvakZwart2007}.
  \subsection{Metropolis algorithm}\label{Metropolis}
  The Metropolis algorithm generates the following Markov chain with stationary density $\pi(x)$
  \begin{eqnarray*}
  X_{i+1}&=&X_{i}+Z_{i+1}\Ii\{U_{i+1}\ge \alpha(X_{i}, X_{i}+Z_{i+1})\}
  \end{eqnarray*}
  for integer $i\ge 0$, where
  \begin{eqnarray*}
  \alpha(x,y)&=& \left\{
                   \begin{array}{ll}
                     \min\{\pi(y)/\pi(x), 1\}, & \mbox{if}\qquad \pi(x)>0, \\
                     1, & \mbox{if}\qquad \pi(x)=0
                   \end{array}
                 \right.
  \end{eqnarray*}
is the acceptance probability to move from $X_i$ to $X_i+Z_{i+1}$ and $X_{i+1}=X_i$ otherwise, \cite{Roberts}. The $Z_i$ has an arbitrary proposal density $q(x)$ and  $U_i$ is a uniformly distributed random variable on $(0,1)$.  Here $\{Z_i\}$ and $\{U_i\}$ are independent sequences of independent, identically distributed random variables, independent of $X_0$. The starting point $X_0$ can be selected arbitrary. The Metropolis algorithm is a special case of the Metropolis-Hastings algorithm with a symmetric proposal density.
\\
Samplings may be compared by  rates of convergence to $\pi(x)$ of the  transition probabilities of a corresponding Markov chain
\[\PP^n(x,A)=\PP\{X_n\in A|X_0=x\},\qquad n\in Z_+,\] to fall  after $n$ steps to a Borel set $A$.
\begin{definition} The Markov chain is called geometrically ergodic if there exists $\rho>1$ such that
\[\lim_{n\to\infty}\rho^n\|\PP\{x,\cdot\}-\pi(x)\|= 0, \]
where $\|\mu\|=\sup_{f: |f|\le 1}|\mu(f)|$ is a total variation norm for a signed measure $\mu$.
\end{definition}
\begin{definition} A Markov chain has polynomial convergence rate $v$ if
\[v=\sup\{\delta: \lim_{n\to\infty}n^{\delta}\|\PP\{x,\cdot\}-\pi(x)\|=0\}.\]
\end{definition}
It is remarkable that the Metropolis Markov chain may have a geometric rate if $\pi(x)$ is a light-tailed distribution, \cite{MengersenTweedie} and, a polynomial rate if $\pi(x)$ is  heavy-tailed, \cite{JarnerRoberts}.
\\
The convergence rate strongly depends on $q(x)$ that can be selected. It is derived for power law that
\begin{equation}\label{4}\pi(x)=\ell(x)/x^{1+r}, \qquad r>0,\end{equation}
the polynomial convergence rate of the Metropolis random walk may be faster ($v=r/\eta$) if \begin{equation}\label{5}q(x)=\ell_q(|x|)/|x|^{1+\eta}, \qquad 0<\eta<2,\end{equation}holds and  slower ($v=r/2$) if $q(x)$ has a finite variance and $0<\lim\inf_{x\to\infty}\ell_q(x)\le \lim\sup_{x\to\infty}\ell_q(x)<\infty$ holds, \cite{JarnerRoberts}. Here, $\ell(x)$ and $\ell_q(x)$ are normalized slowly varying functions.\footnote{This implies that for any $\delta>0$ there exists $K>0$ such that for $y\ge x\ge K$ $\ell(y)/\ell(x)\le (y/x)^{\delta}$.}
Positive constants, functions $\log x$ or $(p_1/p_2(x))^s$, where $p_1(x)$ and $p_2(x)$ are polynomials of the same order, and $s$ is a real number  provide examples of $\ell(x)$.  Functions $(p_1/p_2(x))^s$ and positive constants give examples of $\ell_q(x)$, but not $\log x$.
 \\
 For a Metropolis algorithm there is a simple relation between the geometric ergodicity and the extremal index $\theta$, \cite{Roberts}. Namely, one should check the value of $\eta$ in the limit
 \[\lim_{x\to\infty}(\log \pi(x))'=-\eta\] proposed in \cite{MengersenTweedie}. The Metropolis algorithm is geometrically ergodic if $\eta>0$. If $\eta=0$ then $\theta=0$.
 It was derived that if $\pi(x)$  is a power law then the geometric ergodicity fails and  $\theta=0$ holds. Hence, the mean  first hitting time  of the Metropolis random walk  is infinite in this case.
 \\
 Selecting the proposal distribution $q(x)$ of a Metropolis algorithm, the tail index $\alpha=1+\eta$ is sufficient to find the appropriate polynomial convergence rate of a Metropolis Markov chain to its stationary distribution $\pi(x)$.


\section{Modeling of real data}\label{Model}
The Enron email and DBLP networks  taken from \cite{snapnets} are investigated. We aim first a checking the  power law model (\ref{3}) for these date sets. 
\\
Let $X_1,...,X_n$ be a random sequence of underlying node degrees with the distribution function $F(x)$ and the density $f(x)$,  and $X_{(1)},...,X_{(n)}$ be the corresponding order statistics.
\subsection{Heavy-tail detection}
It follows from the previous section, that it is important to detect the heaviness of the distribution tail and also to estimate the tail index which shows how heavy is the tail. To this end, we evaluate the mean excess function. It is determined by
\[e(u)=\EE(X-u|X>u),\]
and
\[e_n(u)=\sum_{i=1}^n(X_i-u)\Ii\{X_i>u\}/\sum_{i=1}^n\Ii\{X_i>u\}\]
is the sample mean excess function over the threshold $u$. Here, $\Ii\{A\}$ denotes the indicator function of the event $A$.
\\
For heavy-tailed distributions $e(u)$ tends to infinity. Particularly, in case of the Pareto distribution it increases linearly. For light-tailed distributions $e(u)$ tends to zero and it is constant for exponential distribution, \cite{Embr}, \cite{Markovich-2007}.
\\
\begin{figure}[htb]
\centering
\includegraphics[width=0.46\textwidth]{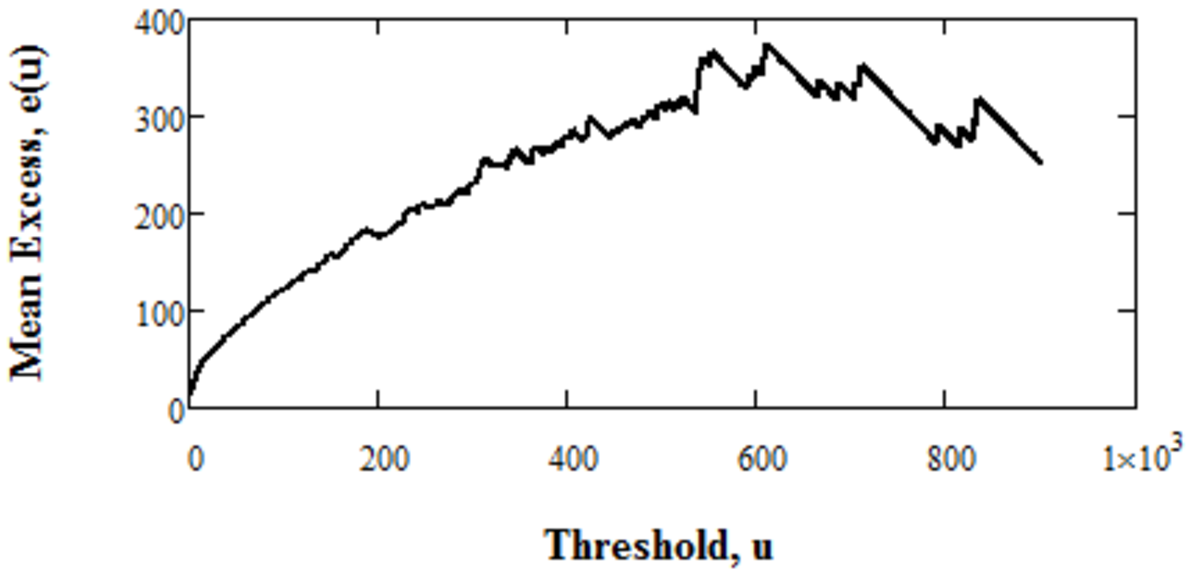}
\includegraphics[width=0.46\textwidth]{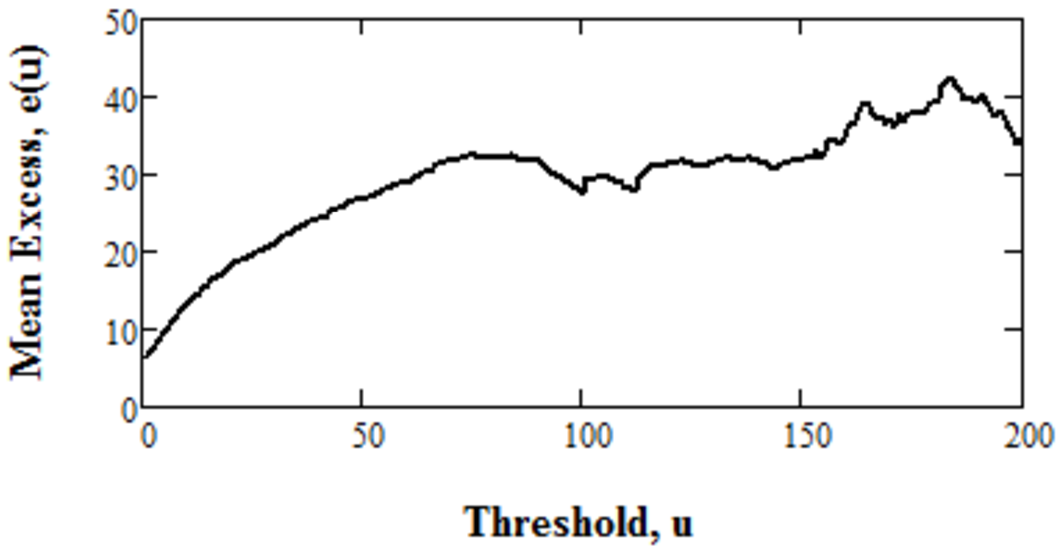}
\caption{Mean excess function of the Enron email data (top) and the DBPL data (bottom).}
\label{fig:4}
\end{figure}
In Fig. \ref{fig:4} one can see plots of the mean excess function
$\{(u,e(u)):X_{(1)}<u<X_{(n)}\}$,  both for Enron and DBPL data sets. We may conclude that both data sets are heavy-tailed. Due to a linearity of the Enron-plot one may suggest that the Pareto  model can be appropriate for the Enron email data. For DBPL data one cannot expect stationarity since the linear curve is changed by a nearly constant line. Hence, we may think that the distribution contains a mixture of an exponential and Pareto distributions. For large $u$ the plots look misleading due to rare observations exceeding such high thresholds.
\\
 The power law density may be determined by
\begin{eqnarray}\label{8}f_{pl}(x)=(X_{(n-k)})^{\alpha}\alpha x^{-\alpha-1},\end{eqnarray} where $\int_{X_{(n-k)}}^{\infty}f_{pl}(x)dx=1$ holds.
\\
The reciprocal of the tail index $\gamma=1/\alpha$  may be estimated by Hill's estimator (\cite{Hill})
\begin{eqnarray}\label{7}\hat{\gamma}^H(n,k)&=&1/k\sum_{i=1}^k \log X_{(n-i+1)}-\log X_{(n-k)}\end{eqnarray}
and by the Ratio estimator. The latter is a generalization of the Hill's one in a sense that instead of  $X_{(n-k)}$ in (\ref{7}) an arbitrary threshold $x_n>0$ is used, \cite{Markovich-2007}. Both estimators may be applied to dependent data such as Markov chains, \cite{Novak}.  Consistency  of Hill's estimator  has been derived in \cite{ResnickStarica} for the $m$-dependent heavy-tailed stationary sequences. We apply also the Moment estimator which is a function
of the Hill's estimator
\[\hat{\gamma}^M_{n,k}=\hat{\gamma}^H(n,k)+1-0.5\left(1-(\hat{\gamma}^H(n,k))^2/S_{n,k})\right)^{-1},\]
 where $S_{n,k}=(1/k)\sum_{i=1}^k\left(\log X_{(n-i+1)}-\log X_{(n-k)}\right)^2$. A survey of other estimators can be found in
  \cite{Markovich-2007} among others.
\\
The number of the largest order statistics $k$ used in all estimators is estimated by a double bootstrap method, \cite{Danielsson}. In Table \ref{tab:0} estimated values of the tail index for the Enron email and the DBLP data are shown. The number of bootstrap re-samples $B=500$ is used.
\\
\begin{table}[!h]
\caption{Estimation of the  tail index $\alpha$  by real data sets averaged over $500$ bootstrap re-samples.}
\begin{tabular}{|c|c|c|c|c|c|}
\hline
Data & Sample & k & Hill & Ratio & Moment
\\
& size &&&&
 \\
\hline
Enron & 36692 & 1659 & 1.337 & 1.2182 & 1.023
\\
DBPL & 425957 & 2589 & 1.028 & 1.277 & 1.657
\\
\hline
\end{tabular}
\label{tab:0}
\end{table}
It is shown that all estimates of $\alpha$  are slightly larger than $1$ but smaller than $2$. This implies the infinite variance of the node degree distribution according to properties of regularly varying distributions (Breiman's theorem), \cite{Embr}, \cite{Markovich-2007}. The fact that not all moments are finite confirms the heavy-tailed distributions of both underlying data sets. Since the tail index is larger for the Enron data, it follows from (\ref{3}) that its distribution has lighter tail than the DBPL data.
%

\begin{figure}[htb]
\centering

\includegraphics[width=0.46\textwidth]{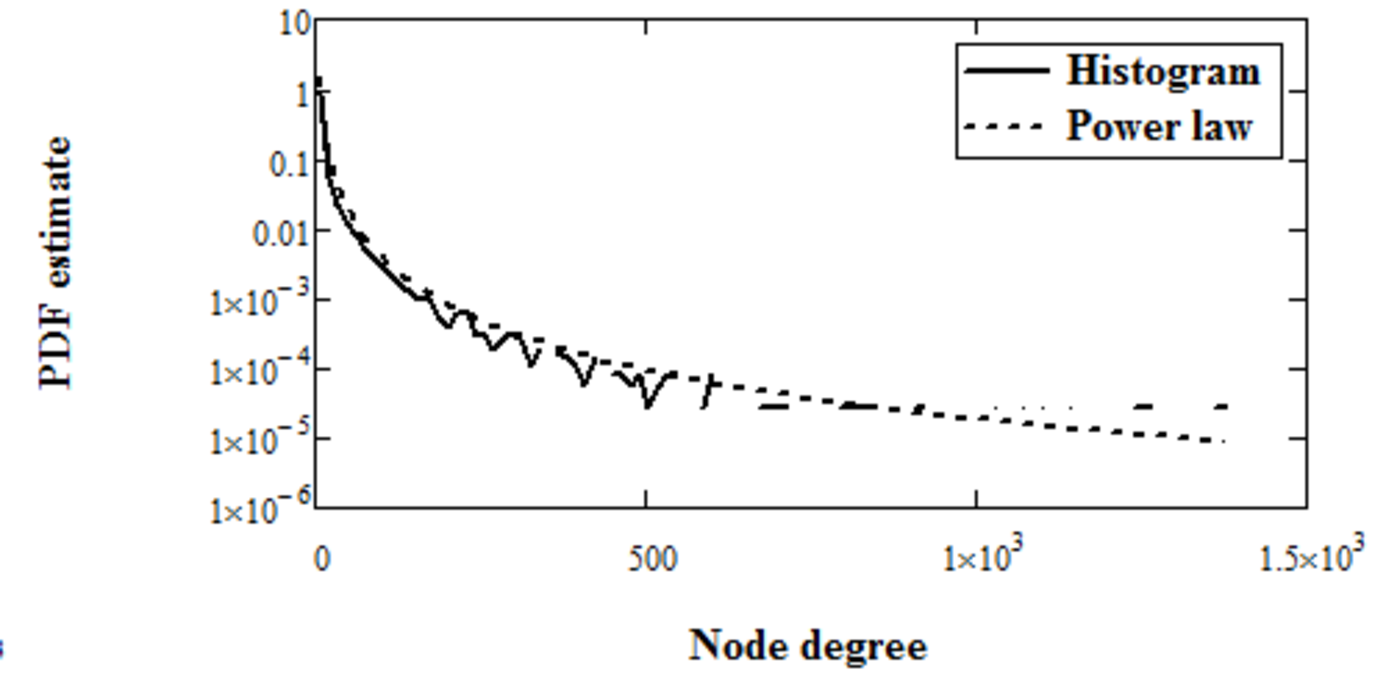}
\includegraphics[width=0.46\textwidth]{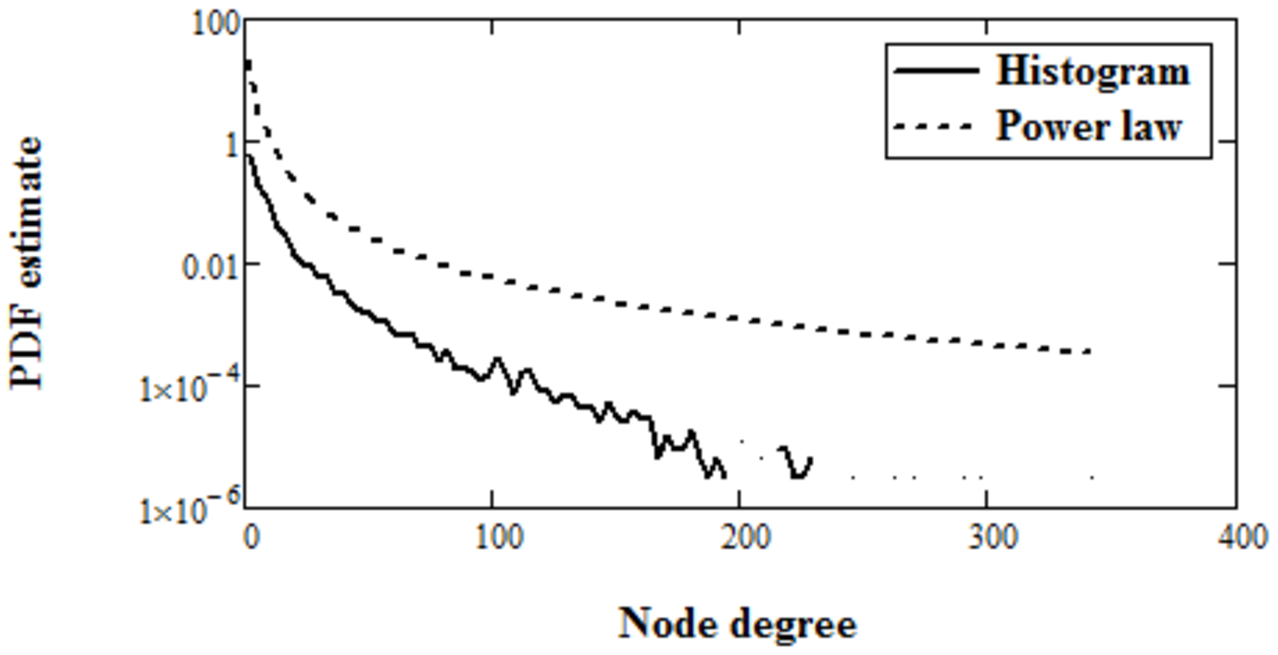}
\caption{The  histogram and the power law density (\ref{8}) for the Enron email data with tail index $\alpha=1.337$ (top) and for the DBPL data with $\alpha=1.277$ (bottom).}
\label{fig:3}
\end{figure}
The comparison of the histogram and the power law density (\ref{8}) for both data sets is shown in Fig. \ref{fig:3}. One may conclude that the power law fits the Enron email data but not the DBPL data. This is in agreement with the mean excess function.
\\
One cannot use the goodness-of-fit tests like
Kolmogorov-Smirnov or von Mises-Smirnov tests to check the hypothesis regarding the power law. The reason  is that the latter tests require samples of independent random variables. However, due to links between nodes, the node degrees are  dependent.
\subsection{Extremal index  estimation}
\begin{figure}[htb]
\centering
\includegraphics[width=0.46\textwidth]{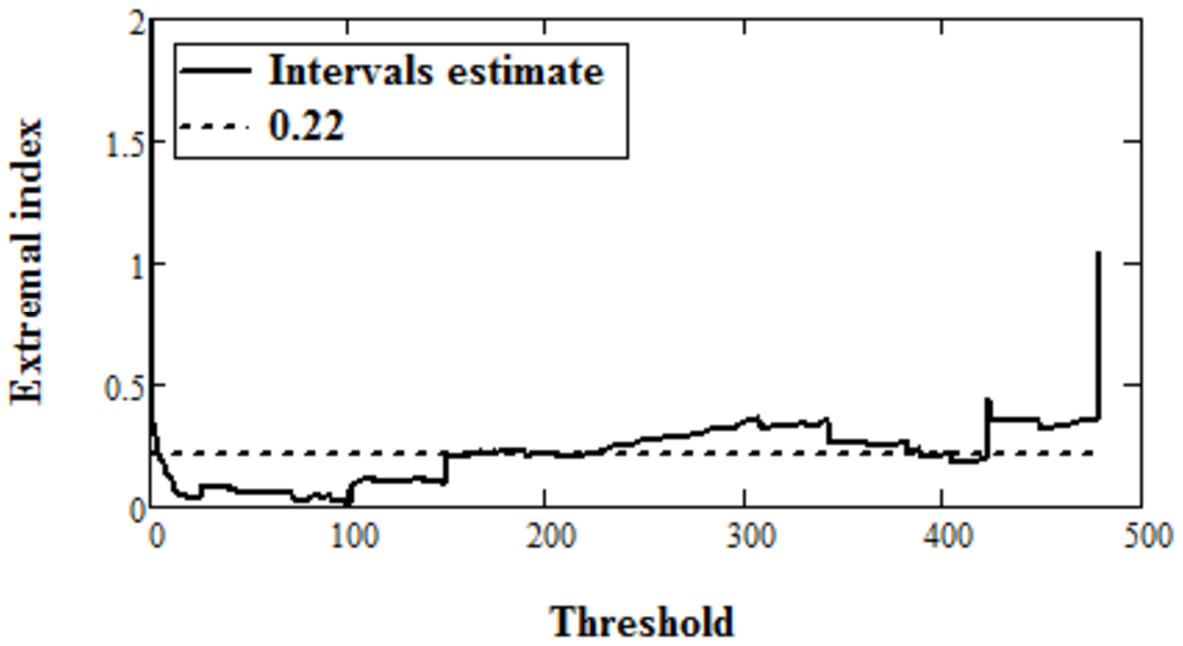}
\includegraphics[width=0.46\textwidth]{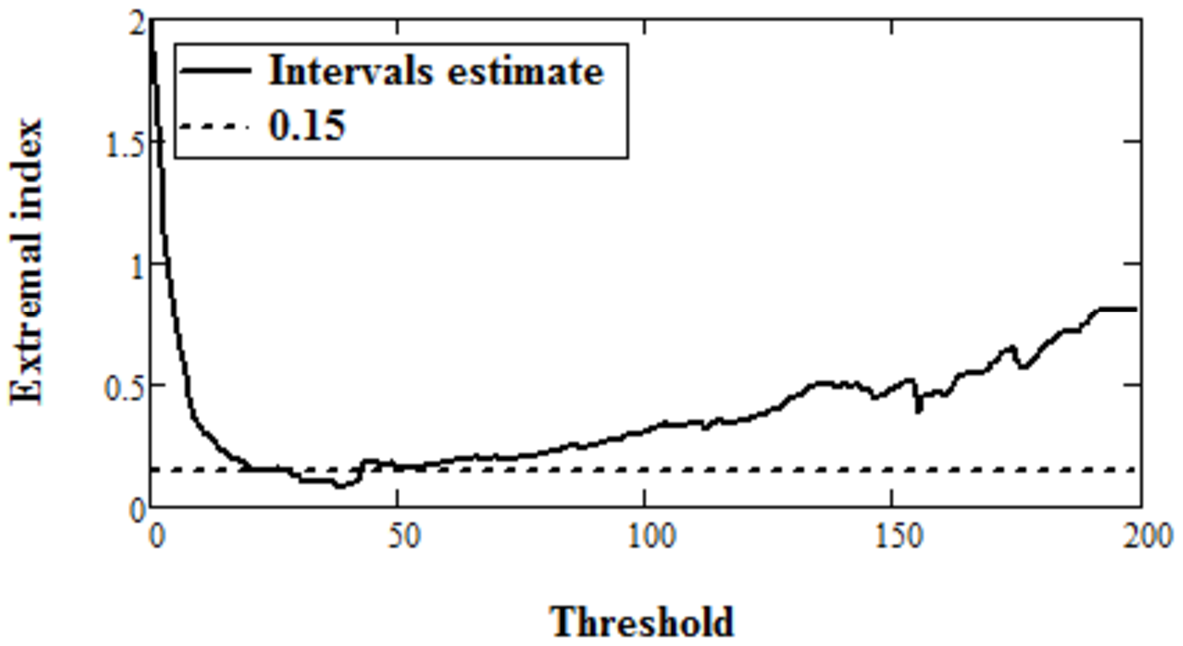}
\caption{Intervals estimate  of $\theta$ for the Enron email data (top) and the DBPL data (bottom): $\theta=0.22$ and  $\theta=0.15$ are selected, respectively, as  values corresponding to stability intervals of the curves.}
\label{fig:8}
\end{figure}

To estimate the extremal index $\theta$, we use the intervals estimator
\begin{eqnarray*}
\hat{\theta}_n(u)&=&\Big\{
\begin{array}{ll}
\min(1,\hat{\theta}_n^1(u)) , \mbox{ if } \max\{T_i : 1 \leq i \leq N -1\} \leq 2,\\
\min(1,\hat{\theta}_n^2(u)) , \mbox{ if } \max\{T_i : 1 \leq i \leq N -1\} > 2,
\end{array}
\end{eqnarray*}
where
\[\hat{\theta}_n^1(u)=\frac{2(\sum_{i=1}^{N-1}T_i)^2}{(N-1)\sum_{i=1}^{N-1}T_i^2},\]
and \[\hat{\theta}_n^2(u)=\frac{2(\sum_{i=1}^{N-1}(T_i-1))^2}{(N-1)\sum_{i=1}^{N-1}(T_i-1)(T_i-2)},\]
proposed by \cite{Ferro}. Here, \[N=\sum_{i=1}^{n} 1(X_i > u)\] is a number of exceedances of $u$ at time epochs $1 \leq S_1 < \ldots < S_N \leq n$ and the interexceedance times are given by $T_i = S_{i+1}- S_{i}$. The intervals estimator does not require the selection of any parameter apart of  $u$ and demonstrates a good accuracy. In contrast,  well-known nonparametric blocks and runs estimators require the size of blocks as an additional parameter to $u$, \cite{Beirlant}.
\\
In Fig. \ref{fig:8} the intervals estimates are shown for both data sets. The appropriate values of $\theta$ are selected  corresponding to a stability interval of the curve $(u, \hat{\theta}_n(u))$. Since $1/\theta$ approximates the mean cluster size of exceedances of the thresholds, we may conclude that the Enron email data contains smaller clusters with $4-5$ nodes on average and the DBPL data $5-6$ nodes on average.
\subsection{First-hitting-time estimation}
Using the obtained estimates $\theta\in\{0.22, 0.15\}$, we evaluate the mean of the first hitting time and its distribution by means of (\ref{9}) and (\ref{11}), respectively.
\begin{figure}[htb]
\centering
\includegraphics[width=0.46\textwidth]{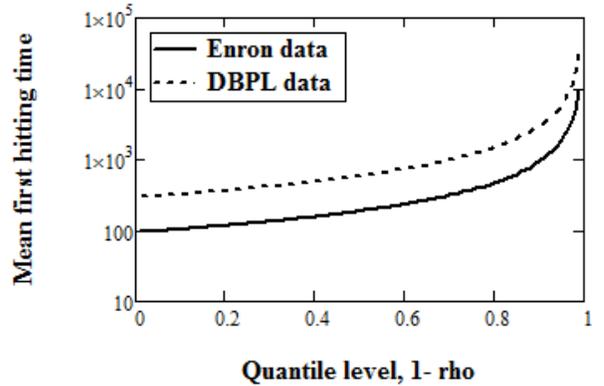}
\caption{The mean first hitting time of the Enron email  and the  DBPL data against the quantile level $1-\rho$: larger $1-\rho$ correspond to larger quantiles $x_{\rho}$ taken as thresholds and longer first hitting times.}
\label{fig:6}
\end{figure}
\begin{figure}[htb]
\centering
\includegraphics[width=0.46\textwidth]{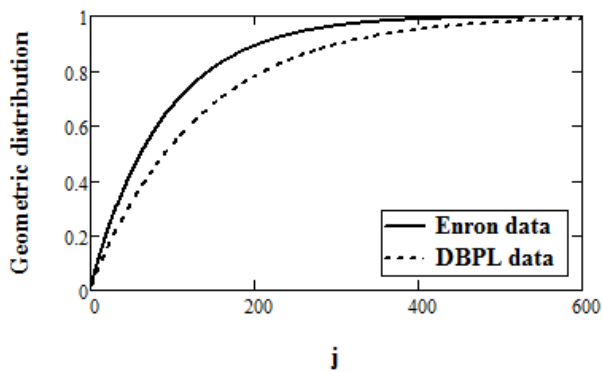}
\caption{The geometric model (\ref{12}) of the normalized distribution $\theta^2P\{T^*(x_{\rho})=j\}$ of the first hitting time of the Enron email  and the  DBPL data corresponding to $95\%$ quantiles as thresholds.}
\label{fig:7}
\end{figure}
In Fig. \ref{fig:6} and \ref{fig:7} we use the following approximations
\[E\left(T^*(x_{\rho_n})\right)\approx 1/(\rho_n\theta^3)\]
 and
\begin{equation}\label{12}\theta^2P\{T^*(x_{\rho_n})=j\}\approx  \rho_n\theta(1-\rho_n\theta)^{j-1},\end{equation}
respectively, that are valid for sufficiently large sample size $n$.
\\
The $95\%$ quantiles $x_{\rho}$ ($\rho=0.05$) of the node degrees were used as thresholds.
\\
Thus, we may conclude that for the DBPL data the mean time required to reach a node with degree larger than $u=x_{\rho}$ is longer than for the Enron email data. This reflects on the distributions, too. The distribution of the first hitting time of the DBPL data has heavier tail than the one of the Enron data.
\subsection{Polynomial convergence rate of the Metropolis random
walk}
Using (\ref{4}) and (\ref{5}) as well as the tail index estimates, one can easily calculate the polynomial convergence rate of the Metropolis random
walk. For example, for $\alpha=1.337$ w.r.t. the Enron data we get from (\ref{4}) that $r=0.337$. To get the largest polynomial rate $v$ we select $\eta$ in (\ref{5}) as small as possible within the interval $(0,2)$. For $\eta=0.01$ we get $v=33.7$. Similarly, for $\alpha=1.028$ w.r.t. the DBPL data we get $r=0.028$ and for the same $\eta$ we obtain  $v=2.8$. This implies, that for the DBPL data whose tail distribution is heavier than the tail of the Enron, we  get the slower polynomial  convergence rate of the Metropolis random walk. Hence, the sampling by means of the Metropolis random walk could be more effective for the Enron email data rather than for the DBPL data.
\section{Conclusions}\label{Conclus}
We  propose to evaluate the optimality of  samplings in complex networks using new measures such as the extremal index, the distribution of the first hitting time and its mean. Considering  real data of social networks we conclude that a heavier tail of the node degree distribution leads to (1) larger node clusters  around the giant nodes, (2) slower convergence rate of the Metropolis random walk, and (3) a longer first hitting time to reach a large node.

\bibliographystyle{plain}
\bibliography{reference}
\end{document}